\documentclass[11pt]{article}

\usepackage[utf8]{inputenc}
\usepackage[T1]{fontenc}

\usepackage[letterpaper,margin=1in]{geometry}
\usepackage{graphicx}
\usepackage{booktabs}
\usepackage{tabularx}
\usepackage{array}
\usepackage{enumitem}
\usepackage{abstract}
\usepackage{titlesec}
\usepackage{xcolor}
\usepackage{caption}
\usepackage{natbib}
\usepackage[hidelinks,colorlinks=true,linkcolor=blue,citecolor=blue]{hyperref}

\definecolor{xnavy}{HTML}{1B2138}
\definecolor{xrule}{HTML}{2E75B6}

\titleformat{\section}
  {\normalfont\large\bfseries\color{xnavy}}{\thesection.}{0.5em}{}
\titleformat{\subsection}
  {\normalfont\normalsize\bfseries\color{xnavy}}{\thesubsection}{0.5em}{}
\titlespacing*{\section}{0pt}{1.4ex plus .2ex}{0.8ex}
\titlespacing*{\subsection}{0pt}{1.1ex plus .2ex}{0.6ex}

\newcolumntype{L}{>{\raggedright\arraybackslash}X}

\title{\vspace{-2.2em}\bfseries\color{xnavy}\LARGE Explainable Optimization:\\[0.15em]
A Call for Interdisciplinary Action}

\author{%
  Nur\c{s}en Ayd{\i}n\thanks{Warwick Business School, University of Warwick. \texttt{nursen.aydin@wbs.ac.uk}} \and
  \c{S}.~\.{I}lker Birbil\thanks{Amsterdam Business School, University of Amsterdam. \texttt{s.i.birbil@uva.nl}} \and
  \.{I}lker K\"u\c{c}\"ukparlak\thanks{Independent researcher, Istanbul. \texttt{ikucukparlak@gmail.com}} \and
  Altu\u{g} Yal\c{c}{\i}nta\c{s}\thanks{Department of Politics and Economics, Ankara University. \texttt{Altug.Yalcintas@politics.ankara.edu.tr}}%
}
\date{\today}

\begin{document}
\maketitle
\begin{abstract}
\noindent
Operations research and management science models support decisions that
affect patients, workers, citizens, public institutions, and organizations.
Decision-makers, such as clinicians approving surgical schedules, planners
allocating disaster relief resources or managers designing workforce
rotations, increasingly require clear and actionable justifications that
bridge the gap between mathematical optimization outputs and the intuitive
reasoning stakeholders need to trust, contest, and implement recommended
decisions. Yet the field has traditionally evaluated optimization models
through computational criteria such as feasibility, optimality, scalability,
and solution time, while treating explanation as a secondary concern.
Mathematical transparency, provided through access to objectives,
constraints, shadow prices, or sensitivity reports, does not automatically
offer the forms of justification that stakeholders need to understand,
trust, contest, or implement optimization-based decisions. This paper calls
for the development of explainable optimization (XOpt) as a distinct
interdisciplinary area that moves beyond algorithmic efficiency and
incorporates behavioral, cognitive, and pragmatic perspectives to address
this explanatory deficit.

\vspace{0.6em}
\end{abstract}

\section{Introduction}
\label{sec:introduction}

Operations research and management science (ORMS) have become indispensable in
modern decision-making, with applications across critical sectors including
logistics, energy systems, healthcare, finance, and public policy. From vehicle
routing algorithms that coordinate global supply chains to mixed-integer
programming models that schedule operating rooms in major hospitals,
optimization techniques now prescribe actions that affect millions of
individuals daily. The field has achieved remarkable success in automating
complex allocation decisions, reducing costs, and improving efficiency metrics
across industries \citep{calma2021operations,petropoulos2024operational}.

Traditionally, ORMS has evaluated models through mathematical and computational
criteria: feasibility, optimality, scalability, robustness, and solution time.
These criteria remain essential, but they are no longer sufficient for deployed
decision-support systems. A solution that is mathematically optimal may still be
rejected, contested, or misused if the people affected by it cannot understand
why it was produced, why plausible alternatives were rejected, or whether the
underlying trade-offs are legitimate. In this sense, optimization faces an
\emph{explanatory deficit}: a gap between what the model can compute and what
stakeholders need in order to trust, contest, and act upon its recommendations.

An explanation, in its most general sense, is a statement that clarifies the
process behind an event or phenomenon (i.e., the \emph{explanandum}). In the
philosophy of science, an explanation is often linked to causation, where an
explanation reveals the hidden causes beneath observable phenomena [\citeauthor{lawson1997econreal}, \citeyear{lawson1997econreal}, pages 21-23; \citeauthor{maki2002dismal}, \citeyear{maki2002dismal}; \citeauthor{reiss2013phil}, \citeyear{reiss2013phil}, pages 101-115].
Explanations, therefore, do not merely restate results but articulate the underlying relations that render them intelligible. However, explanations are not always purely causal; they may also be contrastive, pragmatic, or context-dependent, depending on the question being asked and the audience receiving the answer. Drawing upon \citeauthor{vanfraassen1980scientific}'s [\citeyear{vanfraassen1980scientific}, page 134] view that explanations are context-dependent answers rather than context-free propositions, we argue that this perspective is especially relevant for algorithmic decision systems. An explanation is not merely a technical description of how a model works, but a communicative act that must be made available to stakeholders in a timely and decision-relevant manner: a targeted transfer of information from the computational model or its developer, the \emph{explainer}, to the human operator or affected stakeholder, the \emph{explainee}. The objective of this informational exchange is not merely epistemic clarity, but appropriate trust calibration, ensuring that human reliance on the algorithm is aligned with the system's actual competence, limitations, and context of use. If explanations are delayed, inaccessible, or detached from the moment of decision, stakeholders may lose confidence in both the recommended solution and the optimization model that produced it.

This distinction matters because access to mathematical structure does not
automatically produce understanding, acceptance, or appropriate reliance.
Practitioners and researchers have often assumed that mathematical transparency,
namely the ability to inspect objectives, constraints, assumptions, or
sensitivity reports, is sufficient for practical adoption. Yet as optimization
systems transition from academic demonstrations to deployed solutions affecting clinical decisions, workforce schedules, or emergency response plans, a fundamental gap emerges between computational capability and human acceptability \citep{debock2024explainable}. The gap between knowledge, attitude, and behavior helps explain why factual information alone does not necessarily lead to behavioral action. Specifically, when an algorithm presents a mathematical explanation, confirmation bias may lead stakeholders to place excessive trust in algorithmic outputs when those outputs align with their prior beliefs or expectations. Conversely, when a previous decision informed by an algorithmic explanation
leads to an unfavorable result, negativity bias may produce disproportionate
distrust and algorithm aversion. More broadly, cognitive biases do not affect
all users uniformly; they interact with prior experience, domain expertise,
institutional trust, and the perceived stakes of the decision. As a result, the
provision of mathematically transparent information alone may be insufficient to
support appropriate trust calibration. Explanatory success depends not only on
the correctness of the information provided, but also on how that information is
interpreted, contextualized, and acted upon by different stakeholders \citep{hartisch2026}.

This point is central to \citeauthor{miller2019explanation}'s [\citeyear{miller2019explanation}] argument that effective trust calibration requires drawing upon the ``vast and valuable bodies of research in philosophy, psychology, and cognitive science of how people define, generate,
select, evaluate, and present explanations.'' Miller further argues that the
experts who best understand technical models are not necessarily best placed to
judge whether explanations are useful for intended users, echoing \citeauthor{cooper2004inmates}'s
[\citeyear{cooper2004inmates}] warning that designers often build systems for
themselves rather than for those who must use them. A strong understanding of
the cognitive processes through which humans interact with explanations is
therefore essential. We argue that this observation is equally applicable to the
field of XOpt. ORMS researchers may understand the formulation, algorithm, and
optimality guarantees, but affected stakeholders often need a different form of
justification: one that connects the recommended decision to their
responsibilities, constraints, values, and possible actions.

ORMS has developed sophisticated diagnostic tools for analyzing mathematical
models. Sensitivity analysis, duality theory, and shadow prices provide modelers
with detailed insights into solution structure and parameter dependence. These
methods reveal how optimal solutions respond to changes in input data and, hence,
offer mathematically rigorous characterizations of model behavior. For decades,
these tools have functioned as the primary mechanisms through which optimization
outcomes are interpreted and validated within the field. However, these
diagnostics suffer from a critical limitation: they function primarily as
instruments for modelers rather than explanations for stakeholders. Shadow
prices, for instance, are powerful for analysts fluent in linear programming
theory, but they communicate little to a nurse questioning her shift assignment,
a driver challenging a delivery route, or a patient trying to understand why
surgery has been delayed. The mathematical transparency that satisfies academic
rigor does not necessarily provide the intuitive, actionable justification that
stakeholders require. The distinction between diagnostic power and explanatory
value represents the conceptual fault line which separates traditional ORMS from
the emerging demands of deployed systems \citep{goerigk2023framework}. Recent work on route selection illustrates this shift: rather than explaining the full optimization model, \citet{schild2025} seek concise explanations for why a traffic-aware route differs from an expected traffic-free alternative, thereby moving from model-centric transparency toward decision-centric explanation.

\begin{figure}[t]
  \centering
  \includegraphics[width=\textwidth]{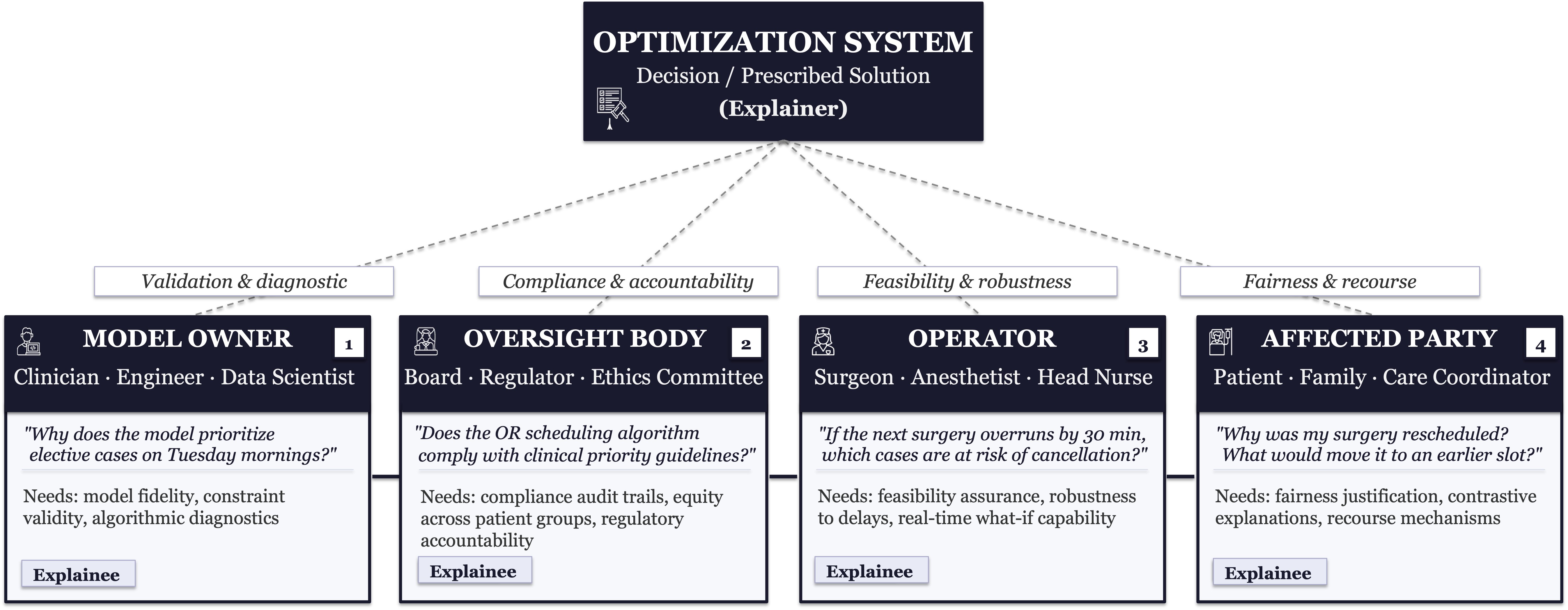}
  \caption{Illustration of the four agent roles in XOpt: model owner, oversight
  body, operator, and affected party.}
  \label{fig:roles}
\end{figure}

We use \emph{stakeholder} as an umbrella term for any individual
or group affected by or involved in an optimization system. We use \emph{agent
role} to denote four broad classes of stakeholders with distinct explanatory
needs: \emph{model owner}, \emph{oversight body}, \emph{operator}, and
\emph{affected party}. As illustrated in Figure~\ref{fig:roles}, these roles
differ not only in their relationship to the optimization system, but also in the
forms of explanation they require. The model owner, typically the designer or
organization, requires explanations that validate model fidelity and identify
improvement opportunities. They need diagnostic tools that reveal whether the
optimization correctly captures business logic and operational constraints. The
operator, such as a dispatcher, nurse, or front-line decision-maker, needs
feasibility and robustness explanations: ``Will this route fail if traffic
increases by 15\%?'' or ``Which constraints are driving this scheduling
decision?'' The affected party, such as a driver, patient, or citizen, seeks
fairness and recourse explanations: ``Why was I assigned this shift?'' or ``What
would change my surgery timing?'' Finally, the oversight body, such as a
regulator or auditor, demands compliance and accountability explanations: ``Does
this solution discriminate against protected groups?'' or ``Can you prove
fairness constraints were enforced?'' In the subsequent discussion, we refer to
the agent providing the explanation as the \emph{explainer} and the agent
receiving it as the \emph{explainee}.

Industrial perspectives on decision-making and recent work on conversational
optimization interfaces confirm that different agents require explanations in
fundamentally different forms---ranging from constraint-activity and trade-off
diagrams to natural-language, counterfactual scenarios or recourse-based actions
\citep{biemans2026explainable,sarpatwar2021explainability, kurtz2026ce4LP, engelhardt2026ce4IP}.
This heterogeneity is not incidental; it reflects structural differences in what
each agent needs from an explanation. Effective explanations must therefore do
more than report why a decision was made. They must clarify why plausible
alternatives were not selected, what trade-offs shaped the outcome, and how the
reasoning should be offered to the relevant explainee. In prescriptive
systems, justification depends not only on the existence of an optimum, but also
on the transparency of the alternatives considered and rejected. This is why a
single explanation format cannot serve all users equally well: the same
optimization result may need to be presented as a diagnostic report, an
operational warning, a fairness justification, or an accountability statement,
depending on who receives it and why.

This perspective implies that an effective explanation cannot be decoupled from
the identification of fundamental user needs. Agent role and expressed intent are
useful indicators, but they may not fully reveal the underlying explanatory need.
For instance, a mathematical structure that satisfies a model owner's need for
system verification will systematically alienate an affected party whose actual
need is actionable recourse. Similarly, the same operator may require different
explanatory structures in different situations: a high-level visual summary during
routine operations, but granular constraint-level information when investigating a
critical failure. Consequently, XOpt must transition toward \emph{need-aware
explanation} frameworks. Such frameworks should preserve the mathematical
integrity of the underlying optimization (diagnostic power), while precisely
addressing the distinct cognitive constraints and functional requirements of the
stakeholder (explanatory value). Rather than relying on static, pre-configured
outputs, XOpt systems should dynamically calibrate explanations to the evolving,
real-time needs of users and decision contexts.

The ORMS community must therefore recognize explanation not as a post-processing
step, but as an intrinsic modeling criterion. Just as constraints encode what is
operationally permissible and objectives encode what is computationally
desirable, explanation requirements encode what is understandable, contestable,
and socially acceptable. Advancing this agenda demands interdisciplinary
collaboration between operations researchers who understand model structure and
algorithmic behavior, behavioral scientists who understand human comprehension
and decision behavior, and domain stakeholders who understand the justification
for the models and data being used to reach certain solutions. Without this triad
of perspectives, optimization risks becoming a technically perfect but socially
impotent discipline.

In our subsequent argumentation, we distinguish interpretability from explainability in XOpt. That is, we refer to interpretability as the transparency of the optimization model, algorithm, or solution structure itself: the reasoning behind the decision can be followed directly from its formulation or design. Explainability, in contrast, refers to the additional means used to help stakeholders understand, evaluate, and question an optimization outcome when the underlying process is too complex to follow. This distinction matters because some optimization systems can be made interpretable by design, while others require explicit explanation mechanisms to make its recommendations understandable, contestable, and usable for a particular stakeholder. For in depth discussion on the concepts of interpretability and explainability in optimization, we refer the interested reader to the recent work of \cite{hartisch2026}.

\section{Call-to-Action: Beyond ``One-Size-Fits-All''}

Recent discussions in ORMS call for complementing traditional computational
metrics, such as solution time and optimality gap, with human-centered metrics
that capture how users comprehend, trust, and act upon algorithmic recommendations
\citep{debock2024explainable}. \citet{hoffman2018metrics} identify
understandability, completeness, usefulness, and trustworthiness as key attributes
of effective explanations and stress that explanations are interactive processes
rather than static outputs. For an explanation to be accepted, it must be tailored
to the user's context, accounting for the triggers that prompt clarification, the
user's existing knowledge, and, most importantly, their underlying goals.

Prior work has classified explanation queries such as ``how,'' ``why,'' ``why
not,'' ``what if,'' to match questions to explanatory methods
\citep{liao2020questioning}, or to group user intentions into classes such as
interpreting a correct, investigating a wrong, or clarifying an ambiguous
prediction \citep{nematov2024aide,adadi2018peeking}. While these frameworks are
valuable, we argue that classifying expressed questions or stated intentions
alone is insufficient. A user lacking domain expertise may articulate an intention
that does not fully align with their actual need \citep{kreminski2024intent}, and
expressed intent may only partially reveal the underlying explanatory need.
Building on the stakeholder-role argument introduced in Section \ref{sec:introduction}, we therefore
call for a need-aware explanation framework in XOpt: a structured process that
(i) diagnoses the underlying psychological need driving an explanation request
through the dimensions described below, and (ii) selects and dynamically
calibrates the explanatory modality to address that need, rather than mapping the
surface form of the question directly to a method.

This need-aware perspective is also consistent with recent empirical evidence on optimization practice. \citet{lawless2025} discuss that optimization model development is not a linear technical pipeline, but an iterative socio-technical workflow involving problem elicitation, data processing, model development, implementation, validation, and deployment. Their study further emphasizes that successful optimization practice is shaped by data, decisions, and sustained dialogue with stakeholders. This supports our argument that explanation in XOpt must be embedded throughout the optimization lifecycle, rather than treated as a final-stage communication layer.

This matters in practice because explanations that miss the underlying need can
actively backfire. Counterfactual explanations, for instance, are widely
recommended for contesting algorithmic decisions, yet their iterative
trial-and-error character imposes substantial cognitive friction on users who are
already frustrated. This, consequently, amplifies dissatisfaction rather than
reducing it. Worse, contrastive explanations based on immutable features
(informing a candidate that they would have been hired had they been younger)
provide no actionable recourse and risk inducing algorithm aversion rather than
acceptance. From the perspective of cognitive load theory, an explanation that is
technically correct can still fail if it overwhelms the working memory of the
explainee \citep{sweller1998cognitive,mayer2003nine}. The same explanation may
therefore be useful for one stakeholder but overwhelming or unhelpful for another,
further illustrating why one-size-fits-all explanatory formats are inadequate.

Diagnosing what a user truly needs requires attention to several structural
dimensions that jointly shape explanatory requirements. System designers should
consider the following:

\begin{itemize}[leftmargin=1.2em,itemsep=0.25em,topsep=0.3em]
  \item[] \textbf{Stakeholder Role:} The agent's position in the system---model
  owner, operator, affected party, or oversight body---provides the baseline
  diagnostic indicator for explanatory requirements, as established in Section \ref{sec:introduction}.
  \item[] \textbf{Cognitive Traits:} Individual cognitive styles affect the
  modality, level of detail, and framing of explanation that a user is likely to
  understand and accept.
  \item[] \textbf{Expertise and Experience:} Domain specialists may require
  probabilistic, granular, or data-grounded explanations, whereas users without
  technical backgrounds may benefit from visual summaries, analogies, or
  plain-language justifications. This suggests the need for layered explanatory
  modalities that adapt to different levels of knowledge and
  information-processing capacity.
  \item[] \textbf{Cultural Structure:} Cultural dynamics, including uncertainty
  avoidance, power distance, and individualistic or collectivist values, may
  shape the required depth, acceptability, and preferred modality of an
  explanation. These forces may also create discrepancies between latent needs
  and expressed intent: users from certain cultural backgrounds may acutely need
  an explanation yet remain conditioned to withhold explicit demand for one.
\end{itemize}

Current research often evaluates isolated explanatory modalities without
examining the interaction between user heterogeneity and explanation
heterogeneity, leaving cognitive impact and practical utility of explanations
insufficiently validated \citep{keane2021ce, longo2024xai}. Relying exclusively on subjective feedback is methodologically insufficient. Such measures should be complemented by behavioral and, where appropriate, unobtrusive indicators of explanation effectiveness, such as decision time, revision behavior, semantic analysis of user responses, or interaction patterns. More broadly, XOpt should move beyond user studies based on hypothetical scenarios toward real-world or high-fidelity simulated settings that test explanations under the complex interplay of user needs, intentions, and emotional states.

Explanations may also take non-causal forms \citep{denboef2025noncausal},
including statistical and pragmatic explanations \citep{vanfraassen1980scientific}.
Counterfactual and contrastive explanations \citep{lewis1973counterfactuals} are
particularly useful in ORMS for comparing a chosen model against competing
alternatives. This clarifies what the outcome would have been under a different
formulation and why certain alternatives were not preferred. This is especially
important when solutions diverge from stakeholder expectations, when no guarantee
of global optimality exists, or when the choice of formulation itself requires
justification. We therefore regard justification as the bridge between
mathematical proof and stakeholder confidence: even a globally optimal solution
requires assurance that the chosen formulation and solution method genuinely
support the decision being recommended.

The ORMS community must therefore design explanation-aware algorithms in which
interpretability is a native constraint rather than a post-hoc patch. Current
practice often separates optimization from explanation: first find a solution,
then attempt to explain it. This modular approach inevitably produces mismatches
between what the optimization computed and what the explanation can articulate.
Embedding explainability during solution generation instead could mean: (i)
adding explicit constraints that bound solution complexity or enforce structural
regularity, such as limiting the number of route segments per driver or using
tree-based decision policies \citep{goerigk2023framework};
(ii) incorporating explanation quality into the objective function, trading
marginal optimality for substantial gains in stakeholder comprehension
\citep{aigner2024framework,forel2023explainable}; and (iii) designing search
procedures that construct solutions through human-traceable logic, such as greedy
heuristics with provable bounds, decompositions aligned with organizational
structure, or explainable heuristic frameworks
\citep{zhou2024,vanstein2025explainable,baczek2025exalt}. These approaches invite a broader
rethinking of fundamental ORMS concepts: the optimality gap could be redefined to
measure distance from both the mathematical optimum and the nearest explainable
solution, and complexity analysis might account for explanation-generation costs
alongside computational costs.

The growth of Explainable Artificial Intelligence (XAI) has raised stakeholder
expectations for algorithmic transparency \citep{danach2025transparent}. While
predictive AI explains what might happen, prescriptive optimization explains what
should be done. This distinction makes explainability even more critical for ORMS.
Yet explainable optimization lags behind XAI in both methodological development
and practical adoption: the field has focused intensely on algorithmic speed and
optimality gaps while treating explainability as secondary. The result is that technically superior
solutions may face resistance, mistrust, or rejection from the humans they are
intended to support. Emerging initiatives, such as the XOpt community
\citep{xopt2025team}, inherently interpretable models, data-driven explainability
frameworks \citep{otto2025coherent}, conversational optimization interfaces
\citep{biemans2026explainable,chen2025optichat}, and recent work on large language models for optimization \citep{wasserkrug2025llm} and reasoning-based model generation \citep{zhang2025decision, xiao2026deepor}, signal the beginning of a more structured effort to close this gap, but these efforts remain fragmented. 

Closing this gap requires rethinking evaluation methodology alongside algorithmic
design. Traditional metrics, such as solution time, optimality gap, and
scalability remain important, but they are not sufficient. XOpt requires
human-centered metrics such as trust calibration, decision-making time with
explanations, error-detection rates when recommendations are flawed, perceived
fairness, contestability, and stakeholder satisfaction across agent roles
\citep{vanstein2025explainable,baczek2025exalt}. In
critical domains such as healthcare and disaster relief, evaluation must
ultimately measure downstream consequences: whether explanations improve
implementation, reduce harmful overrides, support accountability, and increase the
legitimacy of optimization-based decisions. Designing explanation-aware systems
therefore requires expertise in model structure and algorithm design, but also in
how humans perceive and evaluate explanations, and in how concepts such as
fairness, responsibility, and agency are operationalized in real decision
contexts. The call to action is therefore not merely to explain existing models
better, but to develop XOpt as an interdisciplinary field in which explanation
quality is judged against both technical and human-centered criteria.

\section{Case Studies: The Many Faces of Explanation}

This section offers an overview of the earlier conceptual discussion into
real-world examples, demonstrating how challenges and opportunities in XOpt arise
in practice. Rather than providing an exhaustive survey, we focus on three
illustrative cases: (1)~last-mile delivery and workforce fairness, (2)~operating
room scheduling in healthcare, and (3)~disaster relief with supply chain
robustness. Each case highlights different stakeholders, explanation needs, and
moral considerations. In essence, this section shows how the demand for
explainability emerges in each context, where traditional ORMS tools fall short,
and how ideas from XOpt can reshape models, user interfaces, and evaluation
standards. Our appendix complements these case studies by presenting a broader
range of application areas.

\subsection{Last-Mile Delivery and Workforce Fairness}

Last-mile delivery optimization coordinates vast logistics networks through
vehicle routing problems (VRP) and their variants. It employs mixed-integer
linear programming (MILP) as well as advanced optimization approaches such as
column generation within branch-and-price techniques, metaheuristics, and
stochastic or robust methods that can handle uncertain demand and travel times.
The optimization model decides which vehicle serves which customers, in what
sequence, and under which departure times, breaks, capacity limits, and service
constraints, usually with the aim of minimizing cost while satisfying service
targets \citep{toth2014vehicle}.

However, this technical clarity masks an explanatory deficit. Modern VRP solvers
often employ sophisticated techniques that weaken the direct connection between
problem formulation and computed solution. When a driver receives an assigned
route, the path from input data to this specific sequence may involve thousands of
algorithmic iterations: columns are generated, branch-and-bound decisions are
made, neighborhoods are searched, and relaxations are solved. No stakeholder can
mentally trace this algorithmic journey. Recent work on VRP explainability makes
this point explicit: although VRP is widely applied in practice, explainability
for generated routes has remained underdeveloped, even though it is essential for
responsible and interactive route generation. \citet{kikuta2024routeexplainer}
argue that one useful way to explain a route is to explain how each edge
influences the subsequent route, including contrastive questions about why one
edge was selected instead of another. A complementary line of work makes this contrastive perspective even more explicit. \citet{schild2025} study the question that many users of routing systems naturally ask: why is the recommended route different from the one I expected? Rather than explaining the entire routing model, they introduce the notion of a \emph{simple valid explanation}, namely a small set of traffic conditions sufficient to explain why a traffic-aware route differs from the traffic-free alternative. This is particularly relevant for XOpt because it shifts attention from model-centric transparency to decision-centric explanation. In practical terms, a driver or dispatcher often does not need a detailed account of the optimization algorithm; they need a concise and intelligible explanation of why one plausible route was selected instead of another. The approach of \citet{schild2025} therefore illustrates how optimization explanations can be made contrastive, selective, and cognitively manageable for real users.

The multi-stakeholder nature of delivery optimization reveals why
one-size-fits-all explanation fails. Model owners require diagnostic
explanations that help them assess whether the formulation and algorithm behave
as intended: which constraints are active or difficult to satisfy, how time
windows, vehicle capacities, driver-hour limits, and depot restrictions shape the
solution, how sensitive routes are to travel-time estimates, and whether
uncertainty in demand or travel times destabilizes the solution. These technical
explanations satisfy validation needs, but they do not directly address the needs of operators.

Dispatchers managing daily operations need fundamentally different explanations.
They demand feasibility assurance, such as whether a route remains feasible under
traffic disruption; trade-off clarity, such as the cost of prioritizing on-time
delivery over route efficiency; and interactive what-if capability, such as what
happens if a stop is reprioritized or a vehicle breaks down. Most critically, they
need explicit risk articulation, including expected lateness probabilities and
failure modes that algorithmic solutions might ignore.

The most neglected stakeholders could be the drivers and the customers. Drivers may question heavy loads, long routes, and tight time windows: ``Why did I get this assignment rather than that one?'' or ``How was my delivery window determined?'' Customers may want to understand arrival timing, delay handling and service reliability. These explanations require moving beyond optimization diagnostics toward a reasoning where the explainer articulates the value judgments embedded in the model \citep{aigner2024framework,forel2023explainable}. They also
require attention to the broader literature on algorithmic management. The
International Labour Organization emphasizes that algorithmic systems in work
settings can shape working conditions, job quality, control, and fairness
perceptions, making transparency around task assignment and performance evaluation
a matter not only of efficiency but also of worker protection and legitimacy
\citep{baiocco2022algorithmic}.

Let us illustrate. In one of the scenarios, the algorithm proposes Route~B while
the driver anticipated Route~A based on local knowledge and professional
experience. If no explanation is provided, the driver may interpret the assignment
as arbitrary or as evidence that the system disregards their expertise. This can
reduce adoption and acceptability, and may also compromise the driver's sense of
agency and security. A suitable explanation in this instance could be, ``Route~B
was chosen because Route~A is affected by road construction today.'' Such an
explanation does not require exposing the full mathematical model. Instead, it
respects the driver's professional expertise, and provides an actionable reason
for the decision.

This example also illustrates a deeper issue of procedural justice. The
discrepancy between the algorithm's solution of Route~B and the driver's choice of
Route~A is not only a matter of data availability or computational method, it also
concerns epistemic authority. If the optimization model dictates a route without
justification, it risks bypassing the worker's situated knowledge and undermining
their sense of agency. In this context, ``efficiency'' is not a neutral
mathematical objective but a value-laden construct that determines how burdens,
risks, and time pressures are distributed. Providing a reason helps validate the
driver's agency, acknowledge their professional experience, and ensure that the
pursuit of a mathematical optimum does not come at the expense of workplace
fairness or dignity.

The absence of adequate explanation can cascade through the entire logistics
ecosystem. Dispatchers lacking confidence in optimization recommendations may
impose manual overrides, reducing solution quality. Drivers dissatisfied with
route assignments experience lower morale or higher churn, imposing recruitment
and training costs. Here the workforce implications are especially important:
evidence from trucking shows that transparency of algorithmic systems is
associated with justice perceptions, and those justice perceptions mediate
intention to quit \citep{bujold2022opacity}. Missed delivery times can trigger
customer penalties, regulatory exposure for hours-of-service violations, and
reputational harm. Inefficient routes from manual overrides may increase fuel
consumption and emissions, leading to a dual burden of financial and ecological
impacts. Most importantly, unrealistic or poorly explained schedules may create
safety risks, for example when drivers rush to meet tight delivery windows or work
excessive hours despite algorithmic assurances.

\subsection{Healthcare Operating Room Scheduling}

Hospital operating room scheduling demonstrates optimization's critical role in
clinical resource management. The problem coordinates multiple scarce resources,
including operating rooms, surgeons, anaesthetists, nursing teams, and
pre/post-operative beds, to maximize case throughput and on-time surgery starts
while respecting clinical priorities and emergencies. In the ORMS literature, the
problem is typically framed as assigning waiting-list registrations to operating
rooms and shifts subject to procedure duration, specialty compatibility, urgency,
staffing availability, and downstream capacity constraints.

Mathematical formulations may employ again MILP for operating room scheduling and
surgical case assignment, stochastic programming to handle emergency arrivals and
uncertain procedure durations, queueing models to predict bed occupancy, and
robust optimization to manage duration uncertainty. Real-time systems may also
employ rolling-horizon approaches that continuously re-optimize as outcomes
diverge from predictions. More broadly, healthcare planning and scheduling
problems are characterized by many interacting constraints, heterogeneous
resources, and the need to produce usable solutions within short decision times.

From an optimization perspective, delays and cancellations can often be traced to
identifiable bottlenecks: insufficient operating room capacity during peak hours,
limited anaesthetist availability, or insufficient post-operative bed capacity.
Sensitivity analysis can reveal which bottleneck is most critical, for example,
whether adding one anaesthetist would prevent more cancellations than adding one
operating room. These technical insights are valuable for model validation, but
they do not by themselves address the explanatory needs of the people affected by
the schedule.

Model creators require diagnostic explanations about resource bottlenecks and
model behavior. Which constraints bind? How sensitive is on-time surgical
performance to uncertainty in case duration estimates? What is the explicit
trade-off between overtime and cancellations? These questions guide model
refinement, but they do not necessarily help clinical decision-makers 
justify a concrete scheduling change. To make matters more problematic,
\citet{samudra2016scheduling} note a broader issue in the ORMS literature: it is
often unclear whether a study is primarily intended for researchers or
practitioners, and important contextual information is often missing, which makes it harder to assess the managerial relevance of the proposed methods.

Schedulers and clinical leads (operators) demand different
explanations. When cases are moved or delayed, they need accountability: which resource constraints drove the change, and how does the revised schedule respect
prioritization criteria? They require risk articulation: what is the expected
probability of overrun given the current patient mix? They demand what-if
capability: if we added staffing, which cases would benefit? These interactive
explanations require real-time computation beyond batch optimization reports. In
settings such as care unit-constrained scheduling, these explanations are
particularly important because the relevant trade-offs are not only clinical but
also operational, involving utilization, overtime, and downstream congestion.
\citet{fairley2019improving} explicitly note that schedulers may wish to adjust
constraints and re-solve the model, for example by fixing a specific case at a
surgeon's preferred time. In this sense, explanation is closely tied to
contestability and comparison with plausible alternatives, rather than mere
disclosure of an objective value.

Patients and clinical staff are often among the most directly affected
stakeholders, yet they may receive the least explanation. Patients may know that
their surgery has been scheduled, moved, or delayed, but still lack an explanation
of why the change occurred. Did an emergency case take priority? Are bottlenecks
in anaesthetist availability, operating-room capacity, or post-operative beds
preventing earlier scheduling? Do they bear responsibility through excessive
pre-operative fasting, or is this a system limitation? Staff scheduling questions
mirror driver concerns in logistics: ``Why was this shift assigned to me rather
than another anaesthetist?'' and ``How predictable is my shift length given
emergency arrivals?'' Recent work on explainable operating room schedules
reinforces this point directly: \citet{bertolucci2021explaining} argue that an
automated operating room scheduling solution should be explainable to be fully
acceptable in real-world settings, especially when the system fails to produce a
schedule because resources such as beds are insufficient. They further stress that
these explanations should be made available in readable form for inexperienced
users such as medical staff. \citet{mochi2022planning} makes the broader point
that, in healthcare, black-box solutions will be insufficient because both
patients and operators need to know how and why a decision was made.

In this setting, explanation has a psychological as well as technical function. A
probabilistic statement such as ``this case had a high probability of delay'' may
not satisfy a patient or clinician seeking to understand why the disruption
occurred. Explanation supports sense-making, allowing clinicians and patients to
integrate disruptive events---such as a sudden surgery cancellation---into a
coherent and manageable narrative rather than perceiving them as isolated, chaotic
failures. In the absence of such clarity, individuals engage in attribution
processes to assign cause, and often misinterpret systemic resource constraints as
neglect, incompetence, or even malicious intent, whether at the personal or
institutional level. Furthermore, when clinicians, much like other workers, are
required to operate under opaque systems that strip them of a sense of agency over
their labor, this lack of volition frequently manifests as a profound sense of
professional alienation. By providing a transparent causal rationale, an
explanation can restore a degree of perceived control; even when outcomes are
unfavorable, understanding the reason behind a decision can reduce uncertainty and
support appropriate trust.

However, explanation in healthcare scheduling should not stop at mechanical
description. Stakeholders often need to understand not only how a schedule changed,
but why the priority rule behind that change is clinically and ethically
defensible. In this sense, stakeholders need a teleological explanation: an
account that clarifies not only the causes of a scheduling change, but also the
clinical or ethical objective it supports. For instance, if a trauma case occupied
the only available anaesthetist, the patient whose surgery was delayed should be
informed that the decision followed a clinical priority rule rather than
administrative neglect. In this case, explanation functions as justification: it
connects the scheduling outcome to the ethical principles and clinical priorities
that make the decision legitimate. Fairness in clinical scheduling is therefore
not only about mathematically balanced workloads or efficient resource use; it is
also about making sure that every stakeholder (clinical staff, patient, and so on)
can verify that the rules governing their service are rooted in transparent,
justifiable, and just social principles.

The absence of such explanations can create cascading clinical and organizational
costs. Patient dissatisfaction and complaints may damage hospital reputation and
increase legal exposure. Staff may experience burnout when scheduling decisions
appear opaque, leading to turnover in already scarce clinical professions.
Avoidable overtime, whether caused by poor resource allocation or by the inability
to justify decisions to staff, can inflate labor costs. Inefficient resource use
may also follow when staff distrust optimization recommendations and override them.
More broadly, the healthcare explainability literature warns that without methods
for explaining both infeasibility and concrete scheduling decisions, the practical
uptake of AI- and optimization-based scheduling systems will remain limited
\citep{bertolucci2021explaining,mochi2022planning}. In high-stakes environments
such as operating room scheduling, this limitation is not only a technical problem
but also an organizational and ethical one \citep{fairley2019improving}.

\subsection{Disaster Relief and Supply Chain Robustness}

Disaster relief optimization operates under extreme uncertainty and profound
ethical pressure. It spans distinct phases---preparedness, response, recovery, and
mitigation---each with different operational requirements. The challenge is
twofold: deciding where to preposition supplies before disasters whose locations
and magnitudes are uncertain, and determining post-disaster distribution
priorities when demand may vastly exceed available supply. In the humanitarian
logistics literature, these decisions are commonly framed around prepositioning,
facility location, transportation routing, last-mile delivery, and resource
allocation under uncertainty, often across strategic, tactical, and operational
levels. 

From an optimization perspective, the structure appears clear: scenario selection
and probability weights guide robustness decisions, binding access and road
constraints explain distribution bottlenecks, and sensitivity analysis reveals
which uncertainty sources most impact outcomes. Yet the recent humanitarian
logistics research also shows that current AI-enabled decision support remains
concentrated on natural disasters and early-phase response, while recovery,
mitigation, and long-term supply chain resilience remain underexplored
\citep{abdulrashid2026ai}. Persistent challenges include limited explainability of
AI-driven models, siloed system designs, and ethical concerns around fairness and
equity in resource distribution.

Model owners need transparency about foundational modeling and normative
choices. How were scenarios selected and weighted? What trade-off between equity
and efficiency was encoded in the multi-objective formulation? Why were certain
access constraints included or excluded? What assumptions about demand uncertainty
underpin the robust optimization model? These questions matter because, as recent
supply-chain decision-support research argues, explainability is valuable not
merely for technical inspection, but also for enabling informed decision-making
under uncertainty \citep{olan2025enabling}.

Oversight body, including non-governmental organizations (NGOs) and government agencies, require fundamentally different explanations. Why were certain regions prioritized for prepositioning or
post-disaster distribution? What is the expected time to serve affected
survivor communities? If resources shift between regions, what is the impact on coverage,
response time, and equity? Transparency of criteria is essential: are equity
decisions explicit in the model, or do efficiency criteria silently disadvantage
vulnerable populations?

The most critical stakeholders, affected survivor communities, often receive the least
explanation, yet face maximum consequences. Why do some areas receive less
assistance than others? Is this allocation fair? How do vulnerability indicators
influence distribution decisions? Is there a mechanism to appeal or request
assistance if allocation seems inequitable? Conversely, for survivor communities receiving
aid, another question also arises: for how long can support be sustained, and
under what conditions will additional resources be released? These questions show
why explanation in disaster relief must bridge technical optimization, human
justice, and perceived security. The emerging literature emphasizes that equity
and explainability remain comparatively weakly represented objectives in disaster
decision support, even though humanitarian operations require balancing speed,
efficiency, and equity simultaneously under severe uncertainty
\citep{abdulrashid2026ai}.

The risk of secondary disasters makes reserve protection a central feature of
disaster relief optimization. During the response phase, logistics routes may
become unavailable because of undetected infrastructure damage; during recovery,
infectious disease outbreaks or renewed displacement may create additional demand.
From a modeling perspective, holding back critical resources can therefore be a
prudent robustness measure. Yet for affected survivor communities, reserve protection may
be difficult to understand. In the aftermath of trauma and acute scarcity, idle
capacity can be perceived not as precaution, but as delay, neglect, or unequal
treatment. This tension has important behavioral and social dimensions. At the
cognitive level, traumatized individuals may strongly discount future risks
relative to immediate needs. A survivor who lacks shelter, food, medicine, or
clean water may not experience reserve resources as a source of future security,
especially if there is no credible explanation of when and how those resources
will be released. Consequently, when an optimization model holds back aid supplies
to create a buffer against secondary shocks, the affected population might perceive
this idle capacity not as a protective measure but as an inexplicable delay in a
critically needed intervention. At the social level, scarcity can intensify
perceptions of competition between survivor communities. If one region is prioritized over
another, the decision may be interpreted as a direct loss for those who receive
less aid, even when the allocation is based on global coverage, vulnerability, or
robustness criteria. At the institutional level, unexplained withholding or
diversion of resources can be perceived as abandonment by the very authorities or
organizations expected to provide protection. Such perceptions can erode trust and
make coordination more difficult, as survivor communities may resist instructions, hoard
resources, or disengage from official response mechanisms. This example shows that
while optimization may be neutral in its mathematical nature, it cannot remain
neutral on the psychological plane in its social and individual dimensions. To
prevent algorithms from exacerbating existing individual and social
vulnerabilities, the rationale behind these models must be made understandable to
the survivor communities whose recovery depends on them.

As disaster relief logistics illustrates, optimization decisions often carry
distributive meaning. Prioritizing one region, delaying assistance to another, or
holding resources in reserve is not merely an efficiency calculation; it is also
an act of distributive justice. Without explanations of the foundational normative
choices embedded in scenarios, objectives, uncertainty sets, and priority rules,
optimization models risk reinforcing perceptions of arbitrariness or algorithmic
fatalism. To ensure fairness and ethical validity, the relevant assumptions and
decision criteria should be communicated in a form that affected survivor communities,
operational agencies, and oversight bodies can understand and scrutinize.

The failure to explain disaster relief optimization decisions can have severe
humanitarian consequences unmatched in many other optimization domains. Loss of
public trust can undermine the legitimacy of governments and NGOs in crisis
contexts. Social unrest and conflict emerge when survivor communities perceive allocation
as arbitrary or discriminatory. Donor confidence may decline if aid distribution
appears opaque or ineffective. Misdirected aid, delayed response, or inequitable
prioritization can compound suffering. More broadly, we argue that without
interpretable, context-sensitive, and governance-aware decision-support systems,
humanitarian logistics will remain operationally fragmented and difficult to
deploy in practice. In this setting, explanation is not merely a technical add-on;
it is part of accountable and ethically defensible humanitarian decision-making.

\section{From Sensitivity Analysis to XOpt}
Having motivated the need for XOpt and illustrated its many faces through case
studies, we now turn to the methodological toolkit available to the ORMS
community. This section traces a path from legacy analytical instruments that
already carry latent explanatory potential, to methods adapted from explainable
AI, and finally to explanation-native optimization algorithms that embed
interpretability, contestability, or explanation quality into the optimization process itself rather than relying solely on post-hoc explanation add-ons. The goal is not merely to catalogue techniques, but to show how reframing
familiar tools through a stakeholder-centered lens can reorient research
priorities and open new methodological questions for XOpt.

\subsection{Legacy Tools: Power and Limitations}

The ORMS community possesses a rich set of analytical tools that partly serve
explanatory functions, even though they were not originally designed as
stakeholder-facing explanations. These legacy tools emerged primarily from model
validation, sensitivity assessment, and computational analysis. They therefore
provide important diagnostic power, but their explanatory value depends on how
they are translated for different stakeholders.

Duality and shadow prices reveal the marginal value of relaxing constraints. For
example, a shadow price associated with operating-room capacity may indicate how
much the objective would improve if one additional unit of capacity were
available. Such information can be essential for model owners and analysts.
However, this explanation is not automatically meaningful to a surgeon questioning
why a case was moved, or to a patient wondering why their surgery was delayed. The
technical statement must be translated into a stakeholder-relevant justification:
which resource was scarce, how that scarcity affected the decision, and why the
selected schedule was preferred to feasible alternatives.

Sensitivity analysis quantifies how optimal solutions change near reference
parameters. Parametric programming traces optimal objective values across ranges
of input parameters. These tools support analyst-facing ``what-if'' reasoning, for
example by showing how changes in minimum crew-rest requirements affect airline
scheduling costs or how uncertainty in procedure durations affects operating-room
utilization. However, operators and affected parties often ask different
kinds of what-if questions, necessitating a need-aware explanation. A dispatcher
asking ``What happens if I reprioritize this stop?'' needs interactive and
decision-specific explanation, not only sensitivity reports around historical
baseline assumptions.

Scenario generation explores outcome distributions under uncertainty. These
methods are essential for validating robust or stochastic optimization models,
because they reveal how solutions perform across possible futures. Yet presenting
a distribution of possible delays does not necessarily explain why a specific
surgery was scheduled at a specific time, or why a disaster relief model chose to
hold resources in reserve rather than distribute them immediately. Scenario-based
reasoning becomes explanatory only when it connects uncertainty to concrete
decisions, trade-offs, and rejected alternatives.

Comparable taxonomies in explainable AI for OR and industrial
decision-optimization frameworks note that classical sensitivity and duality analysis
fall short of stakeholder-facing explainability: they are technical (requiring
mathematical sophistication), local (explaining behavior near current solutions),
static (producing reports rather than conversations), and diagnostic (validating
models rather than justifying decisions). These tools ask ``Is the model
correct?'' but fail to ask ``Is this decision understandable, contestable, and
legitimate to the relevant stakeholder?''

\subsection{Adapting Methods from Explainable AI}

The growth of XAI offers techniques potentially adaptable to
optimization contexts \citep{danach2025transparent}. However, this transfer must be approached carefully. In the machine learning literature, \citet{rudin2019} warns that post-hoc explanations of black-box models can be unreliable in high-stakes settings and argues that inherently interpretable models should be preferred when comparable performance is achievable. This caution is relevant for XOpt, but the analogy is not exact: optimization models are often mathematically explicit, even when their solutions remain difficult for stakeholders to understand. As \citet{hartisch2026} emphasizes, optimization explanations must address issues that are not central in standard predictive machine learning, including feasibility, active constraints, trade-offs among competing objectives, structural properties of solutions, and changes in recommendations across repeated problem instances.

SHAP (SHapley Additive exPlanations) values, for example, quantify feature
contributions through game-theoretic principles, computing each feature's marginal
contribution across all possible subsets of other features. In optimization contexts, such features are often problem parameters or contextual variables such as resource availability, demand levels, travel-time estimates, cost information, weather conditions, urgency classes, or seasonality. A SHAP-style explanation could therefore be used to assess how these contextual inputs contribute to an optimization output, such as total cost, expected delay, selected route type, or the probability of schedule infeasibility. For example, in vehicle routing, weather conditions, travel-time estimates, customer time windows, and vehicle availability may jointly influence whether a route is selected or whether additional capacity is required. A related but distinct idea is to use SHAP-style attribution to assess the importance of model components, such as constraints, objective terms, or robustness requirements. In this case, the ``players'' in the cooperative game are no longer standard input features, but elements of the optimization model itself. For instance, one may ask how much a time-window constraint, a capacity constraint, or a fairness constraint contributes to the optimal objective value or to the structure of the recommended solution. This interpretation must be handled carefully, because removing, relaxing, or modifying constraints can change the feasible region and may even make some model variants infeasible or operationally meaningless. Thus, adapting SHAP or SHAP-style reasoning to optimization requires explicit attention to feasibility, dependence among constraints, and the distinction between contextual feature attribution and model-component attribution.

Surrogate models offer another route. A simple comprehensible model, such as a
decision tree, may approximate the mapping from input conditions to optimization
outcomes \citep{otto2025coherent}. For example, a surrogate could map demand levels and resource
availability to the number of vehicles required, expected overtime, or the
likelihood of schedule infeasibility. This enables intuitive queries such as ``At
what demand threshold does subcontracting become necessary?'' or ``Which
patient-arrival pattern is likely to force rescheduling?'' Surrogate models can therefore provide
approximate but human-comprehensible descriptions of input--decision relationships when the full optimization process is too complex to communicate directly \citep{goerigk2023framework}.

However, surrogate explanations must be used carefully. A surrogate trained on optimal solutions may misrepresent near-optimal alternatives, obscure feasibility constraints, or suggest decision rules that are not valid for the original optimization problem. For example, a decision tree trained to approximate vehicle-routing outcomes may suggest that two routing patterns are similar in cost, while the underlying optimization model reveals that one violates a time-window or capacity constraint. Thus, surrogate explanations require validation not only in terms of predictive accuracy, but also in terms of feasibility, solution quality, and fidelity to the original optimization model.

A related but distinct approach is data-driven or prototypical explanation. Rather than approximating the optimization process itself, such methods justify a current solution by comparing it with solutions used in similar past instances. For example, \citet{aigner2024framework} propose a framework in which a solution is considered explainable when it resembles favorable solutions implemented under similar circumstances. The explanation therefore rests on historically grounded similarity rather than on a reconstruction of the solver's internal reasoning.

Contrastive and counterfactual explanations are especially promising for
optimization because stakeholders often reason by comparison
\citep{forel2023explainable,kurtz2026ce4LP,engelhardt2026ce4IP}. A driver may ask, ``Why this route
rather than the shorter alternative?'' A production planner may ask, ``Why this
batch size rather than two smaller batches?'' In such cases, the explanation must
compare the selected solution with a plausible alternative and identify the
constraints, costs, or trade-offs that justify the choice. This comparative
structure provides a more intuitive explanation than simply reporting the optimal
objective value.

Implementing contrastive explanations in optimization is itself an optimization
problem. The system must generate alternatives that are meaningfully different from
the recommended solution, feasible under the same constraints, and close enough in
objective value to be relevant \citep{kurtz2026ce4LP}. This requires generating diverse near-optimal
solutions---an inverse problem of optimization where we seek solutions differing
significantly from the recommended one while remaining near-optimal. A closely related contribution is provided by \citet{erwig2024}, who develop a general approach for explaining outcomes of combinatorial optimization algorithms. Their method compares the selected solution with a plausible alternative, or foil, and identifies the specific elements of the solution that account for the difference in performance. This makes it possible to generate concise explanations showing why the chosen solution outperforms the alternative, without requiring the user to inspect the full algorithmic process. Together, work on value-decomposition explanations, explainable data-driven optimization, and coherent local explanations illustrates both the promise and the pitfalls of importing ML-style explainers into optimization without accounting for feasibility, combinatorial structure, and stakeholder-specific cognitive constraints \citep{erwig2024,forel2023explainable,otto2025coherent}.

\subsection{From Explanation Adaptation to Explanation-Native Algorithms}

Adapting XAI methods is valuable but insufficient: XOpt also requires algorithmic
paradigms in which explainability is built into the optimization process itself.
Explanation-native optimization treats interpretability, contestability, and
justification as design criteria, rather than as properties to be recovered after
a solution has already been computed.

One approach is interpretable-by-design, which deliberately trades optimality for interpretability by constraining the solution space to inherently interpretable
structures. Rather than solving a full vehicle routing problem and then
post-processing the solution into explainable form, optimization-by-design
restricts the solution space to self-explanatory routing heuristics and assignment
policies that are easier to communicate. For example, a nearest-neighbor policy is
perfectly self-explanatory: ``The algorithm assigned me to this stop because it is closest
among unserved customers.'' The optimization problem then becomes one of maximizing
performance within this interpretable policy class. The resulting solution may be
slightly suboptimal in aggregate, but each assignment is justified through a
transparent rule.

This idea is formalized in work on inherently interpretable optimization models,
where decision trees encode the decision rule and the optimization problem is
solved under an interpretability constraint \citep{goerigk2023framework}. Such
approaches may sacrifice some optimality, but they offer stronger guarantees of
explainability. Existing work suggests that this trade-off can be valuable in
practice, since deployers may prefer interpretable near-optimal solutions to
mathematically optimal but opaque ones \citep{aigner2024framework}.

A second approach is interactive optimization, which places humans in the
optimization loop and generates explanations through collaborative refinement. This idea is closely related to the interactive multiobjective optimization literature, where decision makers iteratively provide preference information to navigate trade-offs among conflicting objectives \citep{corrente2024}. In an XOpt setting, the system proposes a solution, the stakeholder questions components, and the algorithm responds with targeted explanations while exploring alternatives. For example, a dispatcher may ask, ``Why this driver for this route?'' or a scheduler may ask, ``Could this surgery be moved earlier?'' The system can then explain which constraints become active, what trade-offs are incurred, and why some alternatives
are infeasible or less desirable. This transforms explanations from a static report into a dynamic dialogue with the optimization system. 

Explainable benchmarking for iterative optimization heuristics, interactive multiobjective optimization, and libraries of
explainable algorithmic tools illustrate how analysis of search
behavior, preference information, and alternatives can support human-in-the-loop optimization
\citep{vanstein2025explainable,baczek2025exalt}. When a stakeholder understands
which alternatives are available and why some are preferred, the resulting solution
can gain legitimacy beyond its mathematical optimality. In this sense, interaction
itself becomes part of the explanation process: the system shows relevant
alternatives, explains their trade-offs, and allows stakeholders to participate in
selecting a solution aligned with operational and normative priorities.

Recent work on large language models and reasoning-based optimization systems suggests one possible technological pathway for developing such interactive explanation mechanisms. For example, \citet{zhang2025decision} show how user queries about changes in constraints or parameters can be translated into updated optimization models and accompanying explanations, while \citet{xiao2026deepor} demonstrate the value of explicit multi-step reasoning in optimization modeling. These developments point toward future XOpt systems that do not merely compute optimal solutions, but also expose the assumptions, trade-offs, and modeling decisions behind them in forms accessible to different stakeholders. However, such systems should not be understood as substitutes for human judgment. Rather, they should support stakeholder agency by making optimization logic more contestable, interpretable, and open to scrutiny.

A third approach is constraint-based interpretability, which adds explicit constraints to bound solution complexity or enforce structural regularity. This can be viewed as a form of explainable-by-design optimization: rather than explaining an arbitrary optimal solution after the fact, the model is guided toward solutions with stakeholder-comprehensible regularities. In vehicle routing, for example, a constraint-based description may specify that a route should avoid the city center or use a particular delivery region, without prescribing every arc in the route. Similarly, in workforce scheduling, constraints on the number of shift types, handovers, or assignment patterns can make the solution easier to communicate: ``drivers work within one of three standard shifts'' is more interpretable than a schedule with arbitrary shift boundaries. For crew rostering, constraining solution structure to respect
seniority rules (senior pilots get first choice of preferred pairing
configurations) embeds fairness as an algorithmic primitive rather than a
post-processing filter. Such constraints may reduce optimality, but the interpretability gains can justify the performance cost; in some cases, they may also enable faster algorithms by
restricting the effective search space.

A fourth approach is objective-based explainability, in which explanation quality
is incorporated directly into the objective function, allowing the model to trade
off optimality against comprehensibility. Instead of minimizing operational cost
alone, the model can include a secondary term that penalizes solutions that are
difficult to explain, unstable, fragmented, or misaligned with stakeholder
expectations. For instance, routes with fewer turns at each stop are easier to
communicate to drivers; solutions with fewer binding constraints are easier to
explain; assignments that align with worker preferences require less justification
than counterintuitive assignments. This can be formulated as minimizing operational cost plus a weighted measure of explanation difficulty. The weight attached to this term is not merely a technical parameter; it reflects a value judgment about how much efficiency the organization
is willing to trade for transparency, fairness, and interpretability. This
calibration process can itself become an opportunity to make the value system
embedded in the optimization model explicit.

Together, these approaches invite a broader rethinking of standard ORMS concepts.
The optimality gap could be complemented by an explanation gap: the distance
between the mathematically optimal solution and the nearest solution that is
explainable, implementable, and acceptable to relevant stakeholders. Similarly,
complexity analysis could account not only for solution time, but also for the cost
of generating explanations, exploring alternatives, and supporting interaction. In
this sense, explanation-native algorithms extend classical optimization criteria by
making interpretability, contestability, and justification part of the design space.

\subsection{Evaluation: Beyond Computational Metrics}

Legacy ORMS evaluates optimization solutions primarily through computational
metrics such as solution time, optimality gap, and scalability to larger instances.
These metrics remain essential, but they do not fully capture the objective of
deployed optimization systems: enabling stakeholders to make better decisions and
act on recommendations appropriately \citep{danach2025transparent}.

Human-centered metrics are therefore needed to assess whether explanations achieve
their intended purpose. Trust calibration measures whether stakeholders' confidence
in the solution aligns with its actual quality: overconfidence can lead to
uncritical acceptance, whereas underconfidence can lead to unnecessary overrides
and loss of efficiency gains. Decision-making time quantifies how long stakeholders
require to understand, evaluate, and commit to a recommendation. Error detection
rates measure whether stakeholders can identify when optimization outputs are
suboptimal, infeasible or inappropriate for the decision context
\citep{danach2025transparent}. Other relevant measures include perceived fairness,
contestability, willingness to implement, and satisfaction across different agent
roles.

Recent work on explainable benchmarking, systematic reviews of explainable
decision-making systems, and explainable algorithmic tools highlights the need to
move beyond solver-centric metrics toward user-centered and governance-oriented
criteria \citep{vanstein2025explainable,baczek2025exalt}. Measuring stakeholder
satisfaction requires rigorous experimental design: simply asking whether users are
satisfied is insufficient, because explainability directly influences satisfaction
through comprehensibility (cognitive load), perceived trustworthiness, and fairness.
Evaluation should therefore combine subjective feedback with behavioral indicators,
such as intervention frequency, override patterns, decision time, and the ability to
compare or contest alternatives.

For critical domains including healthcare and disaster relief, evaluation must also
consider downstream consequences. Do explained optimization recommendations improve
patient communication, staff acceptance, or operational throughput? Do transparent
disaster relief allocations improve perceived legitimacy, reduce conflict, or
support better coordination? Do explanations help decision-makers recognize when a
model should not be followed? These outcome measures provide the strongest evidence
that explainability improves decision quality rather than merely improving user
perception.

The research agenda for XOpt therefore requires ORMS to embrace human-centered
evaluation alongside computational performance. A solution that achieves optimal
cost but cannot be understood, justified, or implemented by relevant stakeholders
may fail in practice. Conversely, a slightly suboptimal solution that is
transparent, contestable, and trusted may produce better organizational and social
outcomes. The challenge for XOpt is to make this trade-off explicit, measurable, and
methodologically rigorous.

\section{Conclusion}

In an era characterized by the rapid expansion of algorithmic labor, the necessity
of explainability extends far beyond the superficial pursuit of user comfort,
carrying profound implications for both individuals and society. At the
individual level, this necessity is fundamentally driven by innate psychological
needs, particularly the need for sense-making and perceived control. Without
these, workers may experience professional alienation as they lose agency over
decisions that shape their labor. When optimization systems fail to provide
explanations aligned with these psychological needs, affected individuals may
engage in attribution processes that lead them to misinterpret algorithmic
decisions as neglect, incompetence, or malicious intent on the part of the operator
or model owner. This psychological dimension underscores a critical reality:
providing factual information does not necessarily result in behavioral action or
appropriate trust calibration, largely due to the pervasive influence of cognitive
biases.

Consequently, researchers and practitioners must recognize that the diagnostic
power of an algorithm is not synonymous with its explanatory value. Explainability
cannot be a one-size-fits-all construct; cultural dynamics and different stakeholder
roles necessitate different explanatory approaches. Providing an explanation that is
incongruent with a user's actual needs, or one that exceeds the agent's experience,
domain knowledge, and cognitive capacity, can induce severe cognitive load and
ultimately backfire, resulting in algorithm aversion rather than trust calibration
and informed adoption. A further complicating factor is that human actors may not
always be able to accurately articulate their precise explanatory needs. Therefore,
in both academic research and practical deployment, relying exclusively on
subjective feedback is critically insufficient; evaluations of explanatory efficacy
must also be rigorously grounded in objective, human-centered metrics.

Ultimately, this individual psychological imperative scales to the macro-level.
Explainable optimization represents far more than a technical challenge; it is also
a profound societal and ethical imperative. As ORMS algorithms increasingly
prescribe decisions that fundamentally shape human lives, careers, and safety, the
societal responsibility of researchers expands. The field must evolve from the
narrow pursuit of mathematical objectives to a broader commitment to social and
moral acceptability. The emerging resistance to black-box models is not merely a
technical friction; it is also a demand for justice and a refusal to abandon human
agency in ORMS research.

Beyond technical transparency, explainability also serves as a call for the
democratization of optimization. For too long, the power to define
``optimization'' has been concentrated within a narrow technical community, often
leaving those most affected by algorithmic decisions, including clinicians, workers, and citizens, as passive subjects rather than active participants. By prioritizing
human agency through XOpt, we propose to dismantle the opacity of the ``black-box''
that leads to various forms of biases among cultures, ethnicities, and genders. A
democratized ORMS is one in which the logic of resource allocation, scheduling, and
prioritization is sufficiently accessible to diverse stakeholders. Explainability,
as we see it, is an act of reformation in academic research that returns agency to
the scientific endeavor; it ensures that optimization serves the collective good
rather than just the objective function.

Advances in reasoning foundation models further suggest that explainability may become an integral capability of future optimization systems rather than a separate post-processing step. Recent optimization-focused reasoning models explicitly generate intermediate reasoning traces, exposing the sequence of assumptions, trade-offs, and modeling decisions that lead to a recommendation. Such developments create the possibility of optimization systems that not only prescribe actions but also dynamically generate explanations tailored to the needs of different stakeholders.

ORMS research stands at a philosophical crossroads. Researchers in this field can
continue to produce unexplained or difficult to explain ``black box'' models, or
they can embrace a more mature discipline that integrates explainability from
inception. This choice requires more than ``explainability modules''; it demands
authentic interdisciplinary integration of behavioral science and philosophical
insights. We must recognize that an optimal solution that is not just (or cannot be
justified) is, in a social and economic sense, suboptimal.

We issue this invitation to the INFORMS, EURO, and broader ORMS communities to
recognize explainability as a fundamental research frontier. By treating
explainability as a fundamental design objective, we can do more than improve model
performance; we can better fulfill our duty to the societies that depend on ORMS
researchers. Through XOpt, we can help ensure that optimization serves as a tool for
human flourishing: grounded in transparency, guided by responsibility, and dedicated
to the pursuit of a more just technological future.

\clearpage

\bibliographystyle{plainnat}
\bibliography{references}


\clearpage

\appendix
\section*{Appendix: A Broader Range of Application Areas}
\label{sec:app}
\addcontentsline{toc}{section}{Appendix}

The following tables summarize how the explanatory deficit and the need for
need-aware explanation manifest across a broader range of ORMS application areas.
For each domain, we describe the optimization task, the models and approaches
commonly used, the explanation needs associated with different agent roles (model owners,
oversight body, operator, and affected party), and the consequences of failing to provide adequate explanations.


\newcommand{\apxtablesix}[6]{%
  \par\vspace{2.0em}\noindent\textsc{\large #1}\par\vspace{1.5em}
  \noindent\begin{tabularx}{\textwidth}{>{\bfseries}p{0.40\textwidth} L}
    \toprule
    Short description & #2 \\
    \midrule
    How optimization is used & #3 \\
    \midrule
    Models and approaches & #4 \\
    \midrule
    Explanations demanded & #5 \\
    \midrule
    Consequences of no explanation & #6 \\
    \bottomrule
  \end{tabularx}\par}

\apxtablesix{Last-mile Delivery Vehicle Routing}
{Plan routes for a fleet to deliver orders subject to vehicle capacities, driver-hour limits, and customer time windows.}
{Decide which vehicle serves which customers, in what sequence, with departure times, breaks, and potential split deliveries to minimize cost while meeting service targets.}
{Vehicle Routing Problem variants, typically MILP; column generation/branch-and-price; Lagrangian relaxation; metaheuristics (tabu search, large neighborhood search); stochastic/robust variants for uncertain demand/travel times; dynamic reoptimization with rolling horizon.}
{\emph{(i)~Model owner:} which constraints bind, such as time windows, vehicle capacity, depot resources, and driver time; shadow prices on depot/driver time; sensitivity to travel-time estimates; why certain cuts/columns improve lower bounds; stability under stochastic inputs. 
\emph{(ii)~Oversight body:} compliance with driver-hour regulations, safety standards, service-level commitments, and fairness of workload allocation across drivers or regions. \emph{(iii)~Operator (dispatchers):} why a route is feasible and preferred; trade-offs between cost, workload balance, and on-time delivery; what-if analysis if a stop is reprioritized or a vehicle breaks down; expected lateness risk. \emph{(iv)~Affected party (drivers/customers):} why a driver received a heavy/light load or longer route; why a delivery window is tight; how the estimated arrival time was set; how delays are handled and communicated.}
{Lower operator trust and manual overrides; driver dissatisfaction and churn; missed delivery times and service penalties; inefficient routes and excess fuel; safety risks from unrealistic schedules; regulatory exposure if driver-hour or safety requirements are violated.}

\clearpage
\vspace*{3\baselineskip}

\apxtablesix{Airline Crew Pairing and Rostering}
{Construct legal flight pairings and monthly rosters for pilots and cabin crew under labor, contractual, and aviation rules.}
{Create pairings that cover all flights, then assign crew to rosters while respecting seniority bids, rest times, base constraints, duty-time rules, and fairness considerations with the aim of minimizing cost and disruption.}
{Set covering/partitioning MILP for crew pairing with column generation; rostering as assignment with side constraints; decomposition methods such as Benders and Dantzig--Wolfe; stochastic programming for delays and disruptions; robust optimization for maintenance, weather, and operational uncertainty.}
{\emph{(i)~Model owner:} which legality, rest, base, and coverage constraints are most binding; reduced costs of pairings; why decomposition converges; sensitivity to delay distributions, bid weights, and disruption assumptions. \emph{(ii)~Oversight body:} compliance with aviation regulations, union agreements, duty-time limits, rest requirements, seniority rules, and fairness or non-discrimination requirements. \emph{(iii)~Operator (crew planner):} why a specific pairing/roster was selected; impact of swapping legs or reassigning crew; how cost, fairness, seniority, and disruption resilience were balanced; what-if analysis under delays or aircraft changes. \emph{(iv)~Affected party (crew):} why a bid was not granted; whether rest and duty-time rules are satisfied; fairness across weekends/holidays, night duties, and unpopular pairings; ability to request recourse or understand alternative assignments.}
{Grievances and union disputes; regulatory fines for non-compliance; reduced crew trust and morale; higher sick leave and absenteeism; cascading flight disruptions; reputational harm and compensation costs.}

\clearpage
\vspace*{3\baselineskip}

\apxtablesix{Hospital Operating Room Scheduling and Bed Management}
{Schedule surgeries, assign operating rooms and clinical staff, and plan bed occupancy while accounting for emergencies and uncertainties.}
{Decide case sequencing, operating-room assignment, staffing, and pre/post-op bed allocation to maximize throughput and on-time starts, minimize overtime and cancellations, and respect clinical priorities.}
{MILP for OR scheduling and resource assignment; stochastic programming for emergencies and uncertain case durations; queueing models for bed occupancy; robust optimization for duration uncertainty; rolling-horizon approaches for real-time rescheduling.}
{\emph{(i)~Model owner:} which resource bottlenecks bind such as OR time, beds, anesthetists, or nursing teams; sensitivity to case-duration estimates and emergency-arrival assumptions; trade-offs between overtime, cancellations, and utilization. \emph{(ii)~Oversight body:} compliance with clinical priority rules, patient-safety standards, waiting-list policies, staffing regulations, and fairness across urgency categories or patient groups. \emph{(iii)~Operator (scheduler/clinical lead):} why specific cases were scheduled, moved or delayed; how prioritization criteria were applied; expected risk of overrun; what-if analysis if additional staff, rooms, or beds become available. \emph{(iv)~Affected party (patients/staff):} why a surgery is scheduled/shifted, or canceled; whether the decision follows clinical priority and fairness rules; predictability of shift lengths, breaks, and overtime; opportunities for clarification or recourse.}
{Patient dissatisfaction and complaints; staff burnout and turnover; avoidable overtime and cancelations; inefficient use of operating rooms, beds, and clinical staff; clinical risk from unmanaged bottlenecks; legal or reputational exposure if scheduling decisions appear unsafe, unfair, or poorly justified.}

\clearpage
\vspace*{3\baselineskip}

\apxtablesix{Manufacturing Production Planning and Lot Sizing}
{Plan production quantities and schedules across machines, products, and periods under capacity, setup, and demand constraints.}
{Decide lot sizes, production sequencing, overtime, subcontracting, and inventory levels to meet demand at minimum cost while maintaining service levels and operational stability.}
{MILP for capacitated lot sizing with setup costs/time; sequence-dependent setup scheduling; Lagrangian relaxation; Benders decomposition for multi-stage or network production systems; robust and stochastic programming for demand uncertainty.}
{\emph{(i)~Model owner:} which machine, labour, inventory, and setup constraints bind; shadow prices on capacity constraints; sensitivity to demand forecasts, setup penalties, and service-level requirements; why certain cuts improve solution quality. \emph{(ii)~Oversight body:} compliance with labor agreements, safety rules, customer-service commitments, supplier contracts, and environmental or waste-related constraints. \emph{(iii)~Operator (production planner):} why overtime, subcontracting, or expediting is recommended; trade-offs between inventory, production, and service-levels; impact of machine downtime or rush orders; what-if analysis for forecast changes or capacity disruptions. \emph{(iv)~Affected party (operators/suppliers/customers):} why shifts, workdloads, or delivery schedules changed; rationale for rush orders, delayed production or supplier rescheduling; predictability of workloads and promised delivery dates.}
{Poor operator buy-in and manual overrides; excess inventory or stockouts; missed customer commitments; higher expediting, subconstracting, and overtime costs; supplier dissatisfaction and penalties; reduced workforce trust if schedule changes appear arbitrary or poorly justified.}

\clearpage
\vspace*{3\baselineskip}

\apxtablesix{Retail Inventory Control and Replenishment}
{Set reorder points, order quantities, and stock allocations across stores and distribution centers under uncertain demand.}
{Decide when and how much to order, and how to allocate limited inventory across locations to achieve service-level targets while minimizing holding, storage, waste, and transportation costs.}
{Base-stock policies via stochastic dynamic programming; multi-echelon inventory optimization; newsvendor models; network flow LP for allocation; robust optimization for forecast error; chance-constrainted models for service-level guarantees.}
{\emph{(i)~Model owner:} sensitivity of service levels to lead times, demand variability, and forecast error; which constraints (budget, shelf-space, storage, or transportation) bind; why safety stocks or robust buffers are set. \emph{(ii)~Oversight body:} compliance with allocation policies, consumer-protection requirements, rationing rules, and fairness or non-discrimination expectations across stores, regions, or customer groups. \emph{(iii)~Operator (merchandiser/supply planner):} why a store receives more or fewer units; expected stockout, waste, or markdown risk; what-if analysis if promotions, lead times, or demand forecasts change; assumptions about substitution and product availability. \emph{(iv)~Affected party (store managers/customers):} why items are unavailable, rationed or delayed; expected restock timing; whether allocation decisions are fair across stores, neighborhoods, or customer segments.}
{Lost sales and customer churn; waste, markdowns, and excess inventory; friction between headquarters and stores; perceived inequities in stock allocation; regulatory or reputational scrutiny if rationing or allocation decisions appear biased or poorly justified.}

\clearpage
\vspace*{3\baselineskip}

\apxtablesix{Portfolio Optimization in Asset Management}
{Allocate capital across assets to balance expected return, risk, liquidity, and client mandates under regulatory and investment constraints.}
{Choose portfolio weights subject to risk limits, liquidity requirements, concentration limits, tracking-error bounds, mandate restrictions, and transaction costs; manage rebalancing over time.}
{Mean-variance quadratic programming; conditional value-at-risk (CVaR) minimization via LP; robust optimization against estimation error; scenario optimization; mixed-integer constraints for cardinality/lot sizes, or minimum holdings; multi-period stochastic programming.}
{\emph{(i)~Model owner:} sensitivity to expected returns, covariance estimates, risk models, and transaction-cost assumptions; why regularization or robust uncertainty sets are used; which risk, liquidity, concentration, or mandate constraints bind; stability across scenarios and rebalancing periods. \emph{(ii)~Oversight body:} compliance with investment mandates, regulatory requirements, risk limits, ESG restrictions, fiduciary duties, and auditability of rebalancing decisions. \emph{(iii)~Operator (portfolio manager/compliance):} why portfolio weights changed; contribution of each asset or sector to risk and return; what-if analysis if a constraint is relaxed or market assumptions change; how transaction costs, liquidity, and turnover affected trades. \emph{(iv)~Affected party (clients/investors):} why certain assets were included/excluded; alignment with stated objectives, risk tolerance and ESG preferences; expected drawdown and volatility exposure; clarity on downside protection and long-term suitability.}
{Client mistrust and withdrawals; misalignment with mandates or risk preferences; regulatory and fiduciary concerns; overfitting and hidden concentration risk; excessive turnover and costly whipsaw trading.}

\clearpage
\vspace*{3\baselineskip}

\apxtablesix{Power System Unit Commitment and Economic Dispatch}
{Schedule generation units and dispatch levels to meet electricity demand reliably and at minimum cost under physical, operational, and network constraints.}
{Decide unit startups/shutdowns, ramping, reserve provision, renewable curtailment, and generation dispatch while respecting transmission limits, reliability requirements and renewable integration targets.}
{Mixed-integer unit commitment; DC/AC optimal power flow; stochastic/robust optimization for wind/solar uncertainty; Benders/Lagrangian decomposition; security-constrained optimal power flow.}
{\emph{(i)~Model owner:} which ramping/reserve/transmission/reliability constraints bind; interpretation of dual variables and locational marginal prices; sensitivity to demand, fuel-price, and renewable-generation forecasts; why robustness or security constraints are included. \emph{(ii)~Oversight body:} compliance with market rules, reliability standards, emissions regulations, renewable-integration requirements, and non-discriminatory treatment of generators and consumers. \emph{(iii)~Operator (system/market operator):} why a unit was started, shut down or curtailed; how price components reflect energy, congestion, losses, and reserves; what-if analysis under outages or forecast errors; implications for emissions and system reliability. \emph{(iv)~Affected party (generators/consumers):} why curtailment occurred; why prices spiked; whether dispatch and pricing decisions are fair across locations; assurance that reliability and affordability were considered.}
{Accusations of market manipulation or discriminatory dispatch; regulatory penalties; generator disputes over curtailment or prices; public distrust; blackouts, high prices, emissions impacts, and wider societal costs.}

\clearpage
\vspace*{3\baselineskip}

\apxtablesix{Ride-hailing Matching and Dynamic Pricing}
{Match riders to drivers and set prices in real time to balance supply, demand, service quality, and platform constraints.}
{Assign driver-rider pairs, route drivers, and adjust prices/incentives to reduce waiting times, maintain geographic coverage, and balance rider demand with driver availability.}
{Online bipartite matching and assignment; dynamic pricing via constrained optimization or reinforcement learning with safety constraints; queueing models for waiting times; robust control; multi-objective optimization for fairness, efficiency, and service quality.}
{\emph{(i)~Model owner:} stability and convergence of matching and pricing algorithms; sensitivity to demand forecasts, driver availability, and travel-time estimates; fairness metrics used; constraint activity (driver hours, geographic coverage caps, or surge-price bounds). \emph{(ii)~Oversight body:} compliance with pricing regulations, labor and platform-work rules, consumer-protection requirements, geographic equity expectations, and non-discrimination standards in matching and pricing. \emph{(iii)~Operator (operations team):} why surge pricing was triggered; expected impact on waiting times, driver earnings, and rider acceptance; what-if analysis if price caps, driver incentives, or thresholds change; trade-offs across service regions. \emph{(iv)~Affected party (drivers/riders):} why a price increased; why a driver/rider was matched or not; fairness across neighborhoods, times, and user segments; available recourse or clarification mechanisms.}
{Reputational damage and loss of platform trust; regulatory scrutiny over pricing fairness and worker treatment; driver churn and rider dissatisfaction; inefficient matching, longer waiting times, and reduced service coverage; potential legal action if pricing or matching decisions appear discriminatory.}

\clearpage
\vspace*{3\baselineskip}

\apxtablesix{Facility Location and Logistics Network Design}
{Choose locations for warehouses/hubs and design flows to serve customers under cost, capacity, and service-level requirements.}
{Decide which facilities to open, their capacities, customer-to-facility assignments, and goods flows or routing decisions to minimize total cost while meeting service-level and coverage targets.}
{Capacitated facility location MILP; p-median/p-center models; hub location models; multi-commodity network flow; Benders decomposition; stochastic/robust design for demand, disruption, and lead-time uncertainty.}
{\emph{(i)~Model owner:} which constraints bind (capacity, coverage, service-time); shadow prices on capacity and coverage requirements; sensitivity to demand scenarios, lead-time assumptions, and disruption risks; rationale for robustness. \emph{(ii)~Oversight body:} compliance with planning regulations, environmental requirements, service-access obligations, labor and community-impact commitments, and fairness across regions or customer groups.\emph{(iii)~Operator (network planner/executive):} why a particular site was selected or closed; trade-offs between cost, service level, resilience, and environmental impact; what-if analysis for closing/opening or relocating sites; risk of stranded assets. \emph{(iv)~Affected party (communities/employees/customers):} why a facility opens/closes or relocates nearby; expected effects on jobs, delivery times, service access, traffic, emissions, and local environmental burden}
{Community opposition and implementation delays; poor stakeholder buy-in; misallocation of capital; long-term commitments to suboptimal sites; service inequities across regions; reputational harm.}

\clearpage
\vspace*{3\baselineskip}

\apxtablesix{Disaster Relief Supply Prepositioning and Distribution}
{Plan where to stage relief supplies and how to deliver them after a disaster under uncertain needs and access constraints.}
{Decide prepositioning quantities, depot locations, post-event routing, and prioritization rules to maximize coverage and equity while minimizing time to serve affected populations.}
{Two-stage stochastic programming; robust optimization for scenario uncertainty; equity-aware multi-objective MILP; vehicle routing under damaged networks; chance-constrainted models for coverage and service reliability.}
{\emph{(i)~Model owner:} scenario selection and weighting; assumptions about uncertain demand, infrastructure damage, and resource availability; equity--efficiency trade-offs; binding access, road, depot, or capacity constraints; sensitivity to uncertain demand. \emph{(ii)~Oversight body:} accountability of prioritization criteria; transparency of vulnerability measures; compliance with humanitarian principles, government reporting rules, and equity commitments across regions or population groups. \emph{(iii)~Operator (NGO/government response team):} why certain regions are prioritized; expected time-to-serve; what-if analysis if resources, routes, or access conditions change; rationale for holding reserves or reallocating supplies under secondary-risk scenarios. \emph{(iii)~Affected party (survivor communities):} why some areas receive more, less, earlier, or later assistance; fairness in allocation; how vulnerability indicators influence decisions; opportunities to appeal or request reassessment.}
{Loss of public trust; social unrest; donor withdrawal; misdirected or delayed aid; increased mortality and suffering; legal/political repercussions if allocation decisions appear arbitrary, inequitable, or poorly justified.}

\clearpage
\vspace*{3\baselineskip}

\apxtablesix{Workforce Scheduling in Contact Centers}
{Assign agents to shifts, breaks, and skill groups to meet time-varying call volumes under service-level targets.}
{Decide shift patterns, break times, skill-based routing coverage, overtime, and cross-training plans to minimize labor cost while meeting waiting-time and service-level targets.}
{Integer programming for shift assignment; staffing models with queueing constraints embedded in optimization; multi-skill coverage constraints; stochastic programming for arrival-rate uncertainty; rolling-horizon approaches for real-time schedule adjustment.}
{\emph{(i)~Model owner:} which coverage intervals, skill groups, and staffing constraints bind; sensitivity to arrival-rate estimates; service-level targets, and absenteeism assumptions; why staffing buffers are set; impact of cross-training constraints on cost and service quality. \emph{(ii)~Oversight body:} compliance with labor regulations, contractual rules, break requirements, fairness of shift allocation, service-level reporting, and monitoring of worker impacts.\emph{(iii)~Operator (workforce manager):} why overtime is approved; fairness across nights/weekends; what-if analysis if call volumes change or agents are absent; \emph{(iv)~Affected party (agents/customers):} why agents are assigned to particular shifts, breaks, or skill queues; predictability and fairness of schedules; expected waiting times and service quality for customers.}
{Low morale and attrition; poor schedule adherence; degraded customer service and longer waiting times; increased overtime and operating cost; potential labor disputes if shift allocation appears unfair or opaque; reduced trust in workforce-management systems.}

\end{document}